\documentclass[12pt,titlepage]{elsarticle}

\input xy
\usepackage{amsmath}
\usepackage{amsfonts}
\usepackage{amssymb}

\newcommand{\calH}{{\mathcal H}}

\newcommand{\N}{{\mathbb N}}

\newcommand{\Z}{{\mathbb Z}}
\newcommand{\Q}{{\mathbb Q}}

\newcommand{\C}{{\mathbb C}}

\def\Lyn{{\mathcal Lyn}}

 \def\shuffle{\mathop{_{^{\sqcup\!\sqcup}}}} 

\gdef\stuffle{\;%
  \setlength{\unitlength}{0.0125cm}%
  \begin{picture}(20,10)(220,580) 
  \thinlines 
  \put(220,592){\line( 0,-1){ 10}} 
  \put(220,582){\line( 1, 0){ 20}} 
  \put(240,582){\line( 0, 1){ 10}} 
  \put(230,592){\line( 0,-1){ 10}} 
  \put(225,587){\line( 1, 0){ 10}} 
  \end{picture}\; 
}

\newtheorem{corollary}{Corollary}
\newtheorem{proposition}{Proposition}
\newtheorem{theorem}{Theorem}[section]
\newtheorem{lemma}{Lemma}
\newtheorem{definition}{Definition}

\newtheorem{remark}{Remark}

\newcommand{\Li}{\operatorname{Li}}

\def\deg{\mathrm{deg}}

\newcommand{\poly}[2]{#1 \langle #2 \rangle}

\def\QY{\poly{\Q}{Y}}

\def\QY_0{\Q\left\langle{Y_0}\right\rangle}

\newcommand{\ds}{\displaystyle}

\def\path{\rightsquigarrow}

\gdef\ministuffle{{\scriptstyle \stuffle}}

\def\deg{\mathop\mathrm{deg}\nolimits}

\def\binom#1#2{{#1\choose#2}}

\xyoption{all}

\def\pshuffle#1{{\shuffle}_{#1}} 
\gdef\duffle{\;%
  \setlength{\unitlength}{0.0175cm}%
  \begin{picture}(20,10)(220,580)
  \thinlines
  \put(220,592){\line( 0,-1){ 10}}
  \put(220,582){\line( 1, 0){ 20}}
  \put(240,582){\line( 0, 1){ 10}}
  \put(230,592){\line( 0,-1){ 10}}
  \put(222,587){\line( 1, 0){ 6}}
  \put(225,584){\line( 0, 1){ 6}}
  \put(235,587){\circle*{4}}   
  \end{picture}\;
}

\def\ncp#1#2{#1\langle #2\rangle}
\usepackage{color}
\usepackage{ulem}

\definecolor{amethyst}{rgb}{0.6,0.4,0.8}
  
\definecolor{fuschia}{rgb}{1,0,1} 
 
\definecolor{antiquefuschia}{rgb}{0.57,0.36,0.51} 
 
\definecolor{MyDarkBlue}{rgb}{0,0.08,0.4} 
 
\definecolor{amber}{rgb}{1,0.75,0} 
 
\definecolor{orange}{cmyk}{0,0.6,0.8,0}

\def\mref#1{{\footnotesize ({\ref{#1}})}}

\def\Proof{\medskip\noindent {\it Proof --- \ }}
\def\cqfd{\hfill $\Box$ \bigskip}

\def\scal#1#2{\langle #1 | #2 \rangle}
\def\ncp#1#2{#1\langle #2\rangle}
\def\ncs#1#2{#1\langle \!\langle #2\rangle \!\rangle}
\def\conc{\mathrm{conc}}
\def\span{ \mathrm{span}}

\begingroup
\def\2#1{\ifnum#1<10 0\fi\the#1}
\xdef\isodayandtime{
{\2\day-\2\month-\the\year\space\2{\count0}:%
\2{\count2}}}
\endgroup
\gdef\smuffle{\;%
  \setlength{\unitlength}{0.0125cm}%
  \begin{picture}(20,10)(220,580)
  \thinlines
  \put(220,592){\line( 0,-1){ 10}}
  \put(220,582){\line( 1, 0){ 20}}
  \put(240,582){\line( 0, 1){ 10}}
  \put(230,587){\circle*{5}}
  \end{picture}\;
}

\gdef\shufflemin{\;%
  \setlength{\unitlength}{0.0085cm}%
  \begin{picture}(20,10)(220,580)
  \thinlines
  \put(220,592){\line( 0,-1){ 10}}
  \put(220,582){\line( 1, 0){ 20}}
  \put(240,582){\line( 0, 1){ 10}}
  \put(225,589){\line( 1, 0){ 10}}
  \end{picture}\;
}
\gdef\duffle{\;%
  \setlength{\unitlength}{0.0175cm}%
  \begin{picture}(20,10)(220,580)
  \thinlines
  \put(220,592){\line( 0,-1){ 10}}
  \put(220,582){\line( 1, 0){ 20}}
  \put(240,582){\line( 0, 1){ 10}}
  \put(230,592){\line( 0,-1){ 10}}
  \put(222,587){\line( 1, 0){ 6}}
  \put(225,584){\line( 0, 1){ 6}}
  \put(235,587){\circle*{4}}   
  \end{picture}\;
}
\gdef\huffle{\;%
  \setlength{\unitlength}{0.0175cm}%
  \begin{picture}(20,10)(220,580)
  \thinlines
  \put(220,592){\line( 0,-1){ 10}}
  \put(220,582){\line( 1, 0){ 20}}
  \put(240,582){\line( 0, 1){ 10}}
  \put(225,589){\line( 1, 0){ 10}}
  \put(225,582){\line( 0, 1){ 12}}
  \put(235,582){\line( 0, 1){ 12}}
  \end{picture}\; 
}
\gdef\luffle{\;%
  \setlength{\unitlength}{0.0175cm}%
  \begin{picture}(20,10)(220,580)
  \thinlines
  \put(220,592){\line( 0,-1){ 10}}
  \put(220,582){\line( 1, 0){ 20}}
  \put(240,582){\line( 0, 1){ 10}}
  \put(232,582){\line( 0, 1){ 10}}
  \put(236,587){\circle*{4}}
  \put(223,589){\line( 1, 0){ 5}}
  \put(223,582){\line( 0, 1){ 11}}
  \put(228,582){\line( 0, 1){ 11}}
  \end{picture}\; 
}

\newcommand{\bxi}{\boldsymbol{\xi}}
\newcommand{\brho}{\boldsymbol{\rho}}
\newcommand{\bt}{\mathbf t}
\newcommand{\br}{\mathbf r}
\newcommand{\bS}{\mathbf s}
\newcommand{\calP}{{\mathcal P}}
\def\X{{\mathbb X}}

\newcommand{\pp}{\leqslant}

\usepackage[francais,english]{babel}

\definecolor{MyDarkBlue}{rgb}{0,0.08,0.4}

\usepackage{tikz}
\usetikzlibrary{arrows,decorations.pathmorphing}
\usetikzlibrary{backgrounds,positioning,fit,petri}
\usetikzlibrary{matrix}

\begin{document}


\title{\bf The mechanics of shuffle products\\ and their siblings}

\author{G\'erard H. E.~{\sc Duchamp}${}^\ast$, Jean-Yves~{\sc  Enjalbert}${}^\ast{}$,\\ 
Vincel~{\sc  Hoang Ngoc Minh}${}^\ast{}^\ddag$, Christophe~{\sc  Tollu}${}^\ast$\\[9pt]
\normalsize ${}^\ast$
     LIPN, Institut Galil\'ee - UMR CNRS 7030,\\ 
 93430 Villetaneuse,
France. \\[2mm]
\normalsize ${}^\ddag$
      University Lille 2, 1 Place D\'eliot, 59024 Lille, France. \\[2mm]
}
\begin{abstract}
Nous poursuivons ici le travail commenc\'e dans \cite{DM12} en d\'ecrivant des produits de m\'elanges d'alg\`ebres de fonctions sp\'eciales (issues d'\'equations diff\'erentielles \`a p\^oles simples) de plus en plus grandes. 
Les \'etudier nous conduit \`a d\'efinir une classe de produits de m\'elange, que nous nommons $\varphi$-shuffles. 
Nous \'etudions cette classe d'un point de vue combinatoire, en commen\c{c}ant par \'etendre (sous conditions) le th\'eor\`eme de Radford \`a celle-ci, puis en construisant (toujours sous conditions) sa big\`ebre.  
Nous analysons les conditions des r\'esultats pr\'ecit\'es pour les simplifier en les rendant visible d\`es la d\'efinition du produit de m\'elange.
Nous testons enfin ces conditions sur les produits introduits en d\'ebut d'article.\\[12pt]
{\sl We carry on the investigation initiated in  \cite{DM12} : we describe new shuffle products coming from some special functions
and group them, along with other products encountered in the literature, in a class of products, which we name $\varphi$-shuffle products.
Our paper is dedicated to a study of the latter class, from a combinatorial standpoint. We consider first how to extend Radford's theorem to the products in that class, then how to construct their bi-algebras.
As some conditions are necessary do carry that out, we study them closely and simplify them so that they can be seen directly from the definition of the product.
We eventually test these conditions on the products mentioned above.\\
}
\end{abstract}


\maketitle




\section{Introduction}\label{section-introduction}

As a matter of fact, mathematics (in particular number theory), physics and other sciences provide, for their theories, algebras of special functions indexed by parameters\footnote{The combinatorial supports of these parameters will finally resolve themselves into words.}, with a product, defined at first as a function $X^*\times X^*$ to $\ncp{A}{X}$ and satisfying a simple recurrence of the type
\begin{equation}
 \forall(a,b)\in X^2, \forall (u,v)\in\left(X^*\right)^2,\quad
	au\pshuffle{\varphi} vb=a(u\pshuffle{\varphi} bv)+b(au\pshuffle{\varphi} v)+\varphi(a,b)(u\pshuffle{\varphi} v)\ ,
\end{equation}
the initialization being provided by the fact that $1_{X^*}$ should be a unit.
Of course, we will address the question of the existence of such a product, and will extend it by linearity to $\ncp{A}{X}$.

However, recall that these special functions are indexed by parameters but, unfortunately, sometimes do not exist for some of their values : the prototype of this case is the Riemann zeta function $\displaystyle \zeta(s)=\sum_{n\geqslant 1}\dfrac{1}{n^s}$ for $s=1$.
Nevertheless, if these ``functions'' are seen formally, one can 
in many cases\footnote{That includes in particular all the cases under consideration in our paper}., define a product on the indices which governs the effective product on the functions\footnote{That is the domain of symbolic computation in the vein of Euler and Arbogast\cite{LP,Frechet}.}.

Once the formal identity is obtained, there are many ways to write the divergent quantities as limits of terms which fulfil the same identities (truncated or power series)\footnote{That is the domain of renormalisation and asymptotic analysis initiated by Du Bois-Reymond and Hardy\cite{DBR,Hardy}.}. 

Returning to this family of products, we will use a typology based on examples frequently encountered in the literature as well as new ones that we supply in Section 2. 
\begin{enumerate}
	\item Type I : factor $\varphi$ comes from a product (possibly with zero) between letters (i.e. $X\cup \{0\}$ is a semigroup).
	\item Type II : factor $\varphi$ comes from the deformation of a semigroup product by a bicharacter.
	\item Type III : factor $\varphi$ comes from the deformation of a semigroup product by a colour factor.
	\item Type IV : factor $\varphi$ is the commutative law of an associative algebra (CAA) on $A.X$
	\item Type V : factor $\varphi$ is the law of an associative algebra (AA) on $A.X$
\end{enumerate}
These classes are ordered by the following (strict) inclusion diagram:
\begin{figure}[h]
\centering
\begin{tikzpicture}
  \matrix (m) [matrix of math nodes,row sep=2em,column sep=3em,minimum width=2em]
  {\mbox{I} & \mbox{II} & \mbox{III} & \mbox{V}\\
     &    &  \mbox{IV} &\\
};
\path[right hook->]
    (m-1-1) edge     (m-1-2)
    (m-1-2) edge     (m-1-3)
    (m-1-3) edge     (m-1-4)
    (m-1-2) edge     (m-1-3)
    (m-2-3) edge     (m-1-4);
\end{tikzpicture}
\caption{Hasse diagram of the inclusions between classes.}\label{classes_ord}
\end{figure}
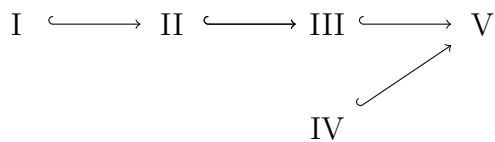
\newpage
We have collected examples from the literature, with the corresponding formulas, in the following table.
\begin{center}
\begin{tabular}{c|c|c|c}\label{tableofexamples}
\hskip-2mm Name\hskip-6mm&\hskip-1mm Formula (recursion)\hskip-6mm&\hskip-1mm $\varphi$\hskip-6mm&\hskip-2mm Type\cr
\hline
\hline
\hskip-2mm Shuffle \cite{ree}\hskip-6mm&\hskip-1mm $au\shuffle bv=a(u\shuffle bv)+b(au\shuffle v)$\hskip-6mm&\hskip-1mm $\varphi\equiv 0$\hskip-6mm&\hskip-1mm I\cr
\hline
\hskip-2mm Stuffle \cite{Hof}\hskip-6mm &\hskip-1mm $x_iu\stuffle x_jv=x_i(u\stuffle x_jv)+x_j(x_iu\stuffle v)$\hskip-6mm&\hskip-1mm $\varphi(x_i,x_j)=x_{i+j}$\hskip-6mm&\hskip-1mmI\cr
&$\phantom{x_iu\stuffle x_jv=}+x_{i+j}(u\stuffle v)$\hskip-6mm&\hskip-1mm\cr
\hline
\hskip-2mm Min-stuffle \cite{C08}\hskip-6mm&\hskip-1mm $x_iu\shufflemin x_jv=x_i(u\shufflemin x_jv)+x_j(x_iu\shufflemin v)$\hskip-6mm&\hskip-1mm $\varphi(x_i,x_j)=- x_{i+j}$\hskip-6mm&\hskip-1mm III\cr
&$\phantom{x_iu\shufflemin x_jv=} - x_{i+j}(u\shufflemin v)$\hskip-6mm&\hskip-1mm\cr
\hline
\hskip-2mm Muffle \cite{EH11}\hskip-6mm &\hskip-1mm $x_iu\smuffle x_jv=x_i(u\smuffle x_jv)+x_j(x_iu\smuffle v)$\hskip-6mm&\hskip-1mm $\varphi(x_i,x_j)=x_{i\times j}$\hskip-6mm&\hskip-1mm I\cr
\hskip-2mm \hskip-6mm &$\phantom{x_iu\smuffle x_jv=}+x_{i\times j}(u\smuffle v)$\hskip-6mm&\hskip-1mm\cr
\hline
\hskip-2mm $q$-shuffle \cite{Bui}\hskip-6mm&\hskip-1mm $x_iu\stuffle_q x_jv=x_i(u\stuffle_q x_jv)+x_j(x_iu\stuffle_q v)$\hskip-6mm&\hskip-1mm $\varphi(x_i,x_j)=qx_{i+j}$\hskip-6mm&\hskip-1mm III\cr
& $\phantom{x_iu\stuffle_q x_jv=}+qx_{i+j}(u\stuffle_q v)$\hskip-6mm&\hskip-1mm\cr
\hline
\hskip-2mm $q$-shuffle$_2$\hskip-6mm&\hskip-1mm $x_iu\stuffle_q x_jv=x_i(u\stuffle_q x_jv)+x_j(x_iu\stuffle_q v)$\hskip-6mm&\hskip-1mm $\varphi(x_i,x_j)=q^{i.j}x_{i+j}$\hskip-6mm&\hskip-1mm II\cr
& $\phantom{x_iu\stuffle_q x_jv=}+q^{i.j}x_{i+j}(u\stuffle_q v)$\hskip-6mm&\hskip-1mm\cr
\hline
\hskip-2mm ${\tt LDIAG}(1,q_s)$ \cite{DTPK}\hskip-6mm&\hskip-1mm\!\!&\!\!\cr
\hskip-2mm (non-crossed,\hskip-6mm&\hskip-1mm $au\shuffle bv=a(u\shuffle bv)+b(au\shuffle v)$\hskip-6mm&\hskip-1mm $\varphi(a,b)=q_s^{|a||b|}(a.b)$\hskip-6mm&\hskip-1mm II\cr
\hskip-2mm non-shifted)&$\phantom{au\shuffle bv=}+q_s^{|a||b|}a.b (u\shuffle v)$\hskip-6mm&\hskip-1mm\cr
\hline
\hskip-2mm $q$-Infiltration \cite{DFLL}\hskip-6mm&\hskip-2mm   $au\uparrow bv=a(u\uparrow bv)+b(au\uparrow v)$\hskip-6mm&\hskip-2mm   $\varphi(a,b)=q\delta_{a,b}a$& III\cr
\hskip-6mm&\hskip-2mm   $\phantom{au\shuffle bv=}+q\delta_{a,b}a(u\uparrow v)$
\hskip-3mm&\hskip-1mm$\phantom{\varphi(a,b)}$ \cr
\hline
\hskip-2mm AC-stuffle\hskip-6mm&\hskip-1mm $au\shuffle_{\varphi} bv=a(u\shuffle_{\varphi}  bv)+b(au\shuffle_{\varphi}  v)$\hskip-6mm&\hskip-1mm $\varphi(a,b)=\varphi(b,a)$ \hskip-6mm&\hskip-1mm IV\cr
\hskip-6mm&\hskip-1mm $\phantom{au\shuffle_{\varphi}  bv=}+\varphi(a,b)(u\shuffle_{\varphi}  v)$
\hskip-6mm&\hskip-1mm $\varphi(\varphi(a,b),c)=\varphi(a,\varphi(b,c))$\hskip-6mm&\hskip-1mm\cr
\hline
\hskip-2mm Semigroup-\hskip-6mm&\hskip-1mm $x_tu\shuffle_{\hskip-1mm\perp} x_sv=x_t(u\shuffle_{\hskip-1mm\perp} x_sv)+x_s(x_tu\shuffle_{\hskip-1mm\perp} v)$\hskip-6mm&\hskip-1mm $\varphi(x_t,x_s)=x_{t\perp s}$\hskip-6mm&\hskip-1mm I\cr
\hskip-2mm stuffle&$\phantom{x_tu\shuffle_{\hskip-1mm\perp} x_sv=}+x_{t\perp s}(u\shuffle_{\hskip-1mm\perp} v)$\hskip-6mm&\hskip-1mm\cr
\hline
\hskip-2mm $\varphi$-shuffle\hskip-6mm&\hskip-1mm $au\shuffle_{\varphi} bv=a(u\shuffle_{\varphi} bv)+b(au\shuffle_{\varphi} v)$\hskip-6mm&\hskip-1mm $\varphi(a,b)$ law of AAU\hskip-6mm&\hskip-1mm V\cr
& $\phantom{au\shuffle_{\varphi} bv=}+\varphi(a,b) (u\shuffle_{\varphi} v)$\hskip-6mm&\hskip-1mm\cr
\hline
\end{tabular}
\end{center}
Of course, the $q$-shuffle is equal to the (classical) shuffle when $q=0$. As for the $q$-infiltration, when $q=1$, one recovers the infiltration product defined in \cite{CFL-4}. 

Many shuffle products arise  in number theory when one studies polylogarithms, harmonic sums and polyz\^etas: it was in order to study all these products that two of us introduced Type IV (see above)  \cite{DM12}. 

On the other hand, in combinatorial physics, one has coproducts with
bi-multiplicative (and noncommutative) perturbation factors (see \cite{GOF18}). 

The structure of the paper is the following:   
in part 2, we complete the first products of \cite{DM12}  with the description of products which come from Hurwitz polyz\^eta functions (the product given in \cite{SLC44} was not valid in all cases) and from generalized Polylerch functions.
We are able to give the complete recursive relation which allows to define all kinds of products; we verify that it implies the existence and uniqueness of this product, which can be extended to $ \ncp{A}{X}$. 
We examine the ``known'' and the ``new'' products in order to determine their classes.
In part 3, we consider how to extend Radford's theorem and we prove that it can be carried over to
the whole class of AC-products (class IV): the Lyndon words constitute a pure transcendence basis of the corresponding commutative algebra, which can moreover be endowed, under additional growth conditions, with a Hopf algebra structure. 

The basis of Lyndon words is the key to effective computations on the algebra of special functions ruled by such products\footnote{The decomposition algorithm (which we shall not decribe in detail) is based on formula (\ref{geneq}) of lemma (\ref{generator}).}.
In part 4, we determine the necessary and sufficient conditions on $\varphi$ so that $\stuffle_\varphi$ belong to the class of AC-products; we give also necessary and sufficient conditions for such a product to be dualizable ({\it i.e.} to be the adjoint of a comultiplication).

\medskip
\textit{Preliminary remark}. It is worth emphasizing at the outset that, although some of the objects/results under review in the present paper have already been defined/proved elsewhere, we include them in our study to lay out as complete a picture as possible and to exemplify the rather 'pedestrian' approach we have adopted. In particular, we have refrained throughout the paper from using more sophisticated algebraic techniques. 

\medskip
\textit{Notation}. In the sequel, $X$ will denote an alphabet, $k$ a $\mathbb{Q}$-algebra, and $A$ a $k$-commutative and associative algebra with unit (a $k$-CAAU).
 

\section{Hurwitz Polyz\^etas and Generalized Polylerch Functions}

We remind the reader of some special functions introduced in \cite{EH11} and complete their study: we prove that they follow a product law which we describe.  

\subsection{Some special functions and their products}

The {\sl Riemann Polyz\^eta} is the function which maps every composition  
$\bS=(s_1,\ldots,s_r)\in (\N_{\geq 1})^r$, to\footnote{The following series converges for $s_1>1$. Under that condition, the definition can be extended by linearity to the module generated by the set of so-called \textit{admissible} composition.} 
\begin{eqnarray}
\zeta(\bS)=\sum_{n_1>\ldots>n_r>0}\dfrac{1}{n_1^{s_1}\ldots n_r^{s_r}}
\end{eqnarray}
 
We now make an observation which, however simple, will appear in different disguises as a building block of many a construction of the paper : There is a (linear) bijection between the module freely generated by (all) compositions and $\ncp{\Q}{Y}$ (where $Y=\{y_k\}_{k\geq 1}$) defined by 
\begin{equation}
\beta_s : (s_1,\ldots,s_r)\mapsto y_{s_1}\ldots y_{s_r}
\end{equation}
So, if $\bS=(s_1,\ldots,s_r)\in (\N_{\geq 1})^r$, $s_1>1$ and $\bS'=(s'_1,\ldots,s'_{r'}), s'_1>1$ are compositions, one knows \cite{EH11} that\footnote{With a slight abuse of language. Stricly speaking, equation \ref{ezs_stuffle} actually reads
$$
\zeta\Big(\beta_s^{-1}\big(\beta_s(\bS)\stuffle\beta_s(\bS')\big)\Big)=\zeta(\bS)\zeta(\bS')\ .
$$
} 
\begin{eqnarray}\label{ezs_stuffle}
\zeta(\bS\stuffle\bS')=\zeta(\bS)\zeta(\bS')
\end{eqnarray}
That function $\zeta$ is well-known and is a special case of the following special functions.

\subsubsection{Coloured Polyz\^etas\\}
The {\sl coloured polyz\^eta} is the function  which, to a composition $\bS=(s_1,\ldots,s_r)$ and a tuple of complex numbers of the same length  $\bxi=(\xi_1,\ldots,\xi_r)$ , associates 
\begin{eqnarray}
\zeta(\bS,\bxi)=\sum_{n_1>\ldots>n_r>0}\dfrac{\xi_1^{n_1}\ldots\xi_r^{n_r}}{n_1^{s_1}\ldots n_r^{s_r}}
\end{eqnarray}
It should be noted that $\zeta(\bS,\bxi)$ appears -- with the notation $\Li_{\bS}(\bxi)$ -- in particule physics \cite{Weinzierl}
.

To describe the product here, we will use two alphabets $Y=\{y_i\}_{i\in\N^*}$, $X=\{x_i\}_{i\in\C^*}$ and $M$ be the (free) submonoid generated by $Y\times X$. One easily checks that\footnote{Throughout the paper $|w|$ stands for the length of the word $w$.} 
$$
M=\{(u,v)\in Y^*\times X^*\ \mid\ |u|=|v|\}
$$ 
As above, to make things rigorous (but slightly more difficult to read), one considers the (linear) bijection defined, on $M$, by 
$$
\beta_c : \left((s_1,\ldots,s_r),(\xi_1,\ldots,\xi_r)\right)\mapsto\left(y_{s_1}\ldots y_{s_r},x_{\xi_1}\ldots x_{\xi_r}\right)\ .
$$
The duffle product is defined as follows. 
\begin{definition} [\cite{DM12}]{\sl (Product of coloured polyz\^etas)} Let $Y=\{y_i\}_{i\in\N^*}$, $X=\{x_i\}_{i\in\C^*}$ and $M$ be as above.\\
The duffle is defined as a bilinear product over $k[M]=\ncp{k}{Y\times X}$ such that 
\begin{eqnarray*}
\forall  w\in M^*,&& w\duffle1_{M^*}=1_{M^*}\duffle w=w,\cr
\forall y_i,y_j\in Y^2, \forall x_k,x_l\in X^2, \forall  u,v\in {M^*}^2,&&
(y_i,x_k).u\duffle (y_j,x_l).v 
=(y_i,x_k)(u\duffle (y_j,x_l)v)\cr
&&+(y_j,x_l)((y_i,x_k)u\duffle v)+(y_{i+j},x_{k\times l})(u\duffle v).
\end{eqnarray*}
\end{definition}
Again, we will show that, under suitable conditions\footnote{Again, rigorously speaking, the left-hand side of the following equation should read
$$
\zeta\left(\beta_c^{-1}\Big(\beta_c(\bS,\bxi)\duffle\beta_c(\bS',\bxi')\Big)\right)\ .
$$
}
\begin{eqnarray}
\zeta\left((\bS,\bxi)\duffle(\bS',\bxi')\right)=\zeta(\bS,\bxi)\zeta(\bS',\bxi')\ .
\end{eqnarray}
\subsubsection{Hurwitz Polyz\^etas}

The {\sl Hurwitz polyz\^eta} is the function  which, to a composition $\bS=(s_1,\ldots,s_r)$ and a tuple of parameters\footnote{\label{ft:parameters} All parameters in the tuple are taken in some subring of $\C$ and none of them is a strictly positive integer.}  of the same length  $\bt=(t_1,\ldots,t_r)$, associates 
\begin{eqnarray}\label{HP}
\zeta(\bS,\bt)=\sum_{n_1>\ldots>n_r>0}\dfrac1{(n_1-t_1)^{s_1}\ldots(n_r-t_r)^{s_r}}\ .
\end{eqnarray}
This series converges if and only if $s_1>1$ (for a ``global'' way to expand (\ref{HP})  as a meromorphic function of $\bS\in\C^r$, see \cite{EH07}). 
To be able to cope with the case $s_1=1$, we have to use the {\sl truncated Hurwitz polyz\^etas} function given by :
\begin{eqnarray}
\forall N\in\N_{>0},\quad\zeta_{_N}(\bS,\bt)=\sum_{N\geqslant n_r>\ldots>n_1>0}\dfrac1{(n_1-t_1)^{s_1}\ldots(n_r-t_r)^{s_r}}
\end{eqnarray}
In order to obtain the product law, we will use here two alphabets $Y=\{y_i\}_{i\in\N_{>0}}$, $Z=\{z_t\}_{t\in{\C\setminus \N_{>0}}}$, the (free) submonoid $N$ generated by $Y\times Z$ 
and, as usual, the bijection
\begin{equation}
\beta_h : \left((s_1,\ldots,s_r),(t_1,\ldots,t_r)\right)\mapsto\left(y_{s_1}\ldots y_{s_r},z_{t_1}\ldots z_{t_r}\right)
\end{equation}
suitably extended by linearity. We have now the following product
\begin{definition}{\sl (Product of Formal Hurwitz Polyz\^etas)} Let $Y=\{y_i\}_{i\in\N^*}$, $Z=\{z_t\}_{t\in k}$ and $N$ be as above.\\
The {huffle}  is defined as a bilinear product over $k[N]=\ncp{k}{Y\times Z}$ such that 
\begin{center}
\begin{tabular}{l}
$\forall  w\in N^*,\qquad\qquad w\huffle 1_{N^*}=1_{N^*}\huffle w=w,$\\[2pt]
$\forall y_i,y_j\in Y^2, \forall z_t,z_{t'}\in Z^2, \forall  u,v\in {N^*}^2,$\\[2pt]
$\displaystyle
\begin{array}{rcl}
t= t'&\Rightarrow&
(y_i,z_t)u\huffle (y_j,z_{t})v\\
&&=(y_i,z_t)(u\huffle (y_j,z_{t})v)+(y_j,z_{t})((y_i,z_t)u\huffle v)\\
&&+(y_{i+j},z_{t})(u\huffle v)
\end{array}
$\\[2pt]
$\displaystyle\begin{array}{rcl}
t\not= t'&\Rightarrow&
(y_i,z_t).u\huffle (y_j,z_{t'}).v\\
&&= 
(y_i,z_t).\left(u\huffle (y_j,z_{t'}).v\right)
+(y_j,z_{t'}).\left((y_i,z_t).u\huffle v\right)\\
&&+\displaystyle\sum_{n=0}^{i-1}{j-1+n\choose j-1} 
\frac{(-1)^n}{(t-t')^{j+n}}\,(y_{i-n},z_t).\left(u\huffle v\right)\\
&&+\displaystyle\sum_{n=0}^{j-1} {i-1+n\choose i-1} 
\frac{(-1)^n}{(t'-t)^{i+n}}\,(y_{j-n},z_{t'}).\left(u\huffle v\right)\ .
\end{array}
$
\end{tabular}
\end{center}
\end{definition}
We also will show that\footnote{Again, rigorously speaking, the left-hand side of the following equation should read
$$
\zeta_N\left(\beta_h^{-1}\Big(\beta_h(\bS,\bt)\huffle\beta_h(\bS',\bt')\Big)\right)\ .
$$
} for all integer $N$
\begin{eqnarray}
\zeta_N\left((\bS,\bt)\huffle(\bS',\bt')\right)=\zeta_N(\bS,\bt)\zeta_N(\bS',\bt')\ .
\end{eqnarray}
\begin{remark}
The functions we call 'Hurwitz polyz\^etas', a term coined in the last century (see for example \cite{SLC44}).
 must not be confused with the monocenter polyz\^etas, defined only for a composition $\bS$ and a parameter $t$ by
\begin{eqnarray}
\zeta(\bS,t)&=&\sum_{n_1>\ldots>n_r>0}\dfrac1{(n_1-t)^{s_1}\ldots(n_r-t)^{s_r}},
\end{eqnarray}
which follow a much simpler rule, namely the stuffle product on the compositions. 
\end{remark}
\subsubsection{Generalized Polyl\^erch functions}
The {\sl generalized Polyl\^erch function} is the function  which maps a composition $\bS=(s_1,\ldots,s_r)$, a tuple $\xi=(\xi_1,\ldots,\xi_r)$ of complex numbers, and a tuple $\bt=(t_1,\ldots,t_r)$ of parameters$^{\ref{ft:parameters}}$, all three of the same length, to 
\begin{eqnarray}
\zeta(\bS,\bt,\bxi)=\sum_{n_1>\ldots>n_r>0}\dfrac{\xi_1^{n_1}\ldots \xi_r^{n_r}}{(n_1-t_1)^{s_1}\ldots(n_r-t_r)^{s_r}} .
\end{eqnarray}
Here, we will need three alphabets $Y=\{y_i\}_{i\in\N^*}$, $X=\{x_i\}_{i\in\C^*}$,  $Z=\{z_t\}_{t\in k}$ and the (free) submonoid $T$ generated by $Y\times Z\times X$. The bijection
\begin{equation}
\beta_l : \left((s_1,\ldots,s_r),(t_1,\ldots,t_r),(\xi_1,\ldots,\xi_r)\right)\mapsto\left(y_{s_1}\ldots y_{s_r},z_{t_1}\ldots z_{t_r},x_{\xi_1}\ldots x_{\xi_r}\right)
\end{equation}
still extended by linearity. The product $\luffle$ is given by the following definition:
\begin{definition}{\sl Product of Generalized Lerch functions}\\ 
Let $Y=\{y_i\}_{i\in\N^*}$, $X=\{x_i\}_{i\in\C^*}$,  $Z=\{z_t\}_{t\in k}$ 
and $T$ be the (free) submonoid generated by $Y\times Z\times X$.\\
The {luffle} is defined as the bilinear product over $k[T]=\ncp{k}{Y\times Z\times X}$ satisfying the following recursive relation :
\begin{center}
\begin{tabular}{l}
$\forall  w\in A^*,\qquad\qquad w\luffle 1_{A^*}=1_{A^*}\luffle w=w,$\\[2pt] 
$\forall (y_i,y_j)\in Y^2, \forall (z_t,z_{t'})\in Z^2, \forall (x_k,x_l)\in X^2, \forall  (u,v)\in {A^*}^2,$\\[2pt]
$\displaystyle
\begin{array}{rcl}
t= t'&\Rightarrow&
(y_i,z_t,x_k).u\luffle (y_j,z_{t},x_l).v\\
&&= 
(y_i,z_t,x_k).\left(u\luffle (y_j,z_{t}).v\right)
+(y_j,z_{t},x_l).\left((y_i,z_t).u\luffle v\right)\\
&&+(y_{i+j},z_t,x_{k\times l}).\left(u\luffle v\right)\\
\end{array}$\\[2pt]
$\displaystyle\begin{array}{rcl}
t\not= t'&\Rightarrow&
(y_i,z_t,x_k).u\luffle (y_j,z_{t'},x_l).v\\
&&= 
(y_i,z_t,x_k).\left(u\luffle (y_j,z_{t'}).v\right)
+(y_j,z_{t'},x_l).\left((y_i,z_t).u\luffle v\right)\\
&&+\displaystyle\sum_{n=0}^{i-1}{j-1+n\choose j-1} 
\frac{(-1)^n}{(t-t')^{j+n}}\,(y_{i-n},z_t,x_{k\times l}).\left(u\luffle v\right)\\
&&+\displaystyle\sum_{n=0}^{j-1} {i-1+n\choose i-1} 
\frac{(-1)^n}{(t'-t)^{i+n}}\,(y_{j-n},z_{t'},x_{k\times l}).\left(u\luffle v\right)
\end{array}$
\end{tabular}
\end{center}

\end{definition}
We also show\footnote{Again, rigorously speaking, the left-hand side of equation \ref{seventeen} should read
\begin{equation}
\zeta\left(\beta_l^{-1}\Big(\beta_l(\bS,\bt,\bxi)\luffle\beta_l(\bS',\bt',\bxi')\Big)
\right)\ .
\end{equation}}
\begin{eqnarray}\label{seventeen}
\zeta\left((\bS,\bt,\bxi)\luffle(\bS',\bt',\bxi')\right)=\zeta(\bS,\bt,\bxi)\zeta(\bS',\bt',\bxi')\ .
\end{eqnarray}

\subsection{General framework of study}

Other products from table \ref{tableofexamples} belong to the same family as the examples examined so far, and so pertain to the same kind of approach. As we aim to offer as comprehensive a framework as possible, we now concentrate on the most general class of $\varphi$-products, \textit{i.e.} class V, which emerges from definition \mref{Defstarproduct} below. 
We will use a unitary ring as the ground set of scalars (and not a field as it would be expected in combinatorics) because some applications require to work with rings of (analytic or arithmetic) functions. 
\begin{proposition}\label{prop:Starproduct}
Let $A$ be a unitary commutative ring, $X$ be an alphabet and $\varphi: X\times X\rightarrow \ncp{A}{X}$ is an arbitrary mapping.
Then there exists a unique mapping $\star: X^*\times X^*\rightarrow \ncp{A}{X}$ satisfying the conditions :
\begin{eqnarray}
(R)
\begin{cases} 
\text{for any} & w\in X^*,\ 1_{X^*}\star w=w\star  1_{X^*}=w,\cr 
\text{for any} & a, b\in X \text{and } u,v\in X^*,\cr 
 & au\star bv 
=a(u\,\star\,bv)+b(au\star v)+\varphi(a,b)(u\star v).
\end{cases}
\end{eqnarray} 
\end{proposition}
\Proof
By recurrence over $n=|u|+|v|$. 
\cqfd
\begin{definition}\label{Defstarproduct}
With the notations of Proposition \ref{prop:Starproduct}, the unique mapping from $X\times X$ to $\ncp{A}{X}$ satisfying conditions (R) will be noted $\pshuffle{\varphi}$ and will be called $\varphi$-shuffle product.
\end{definition}
\smallskip
From now on, we suppose that $\varphi$ takes its values in $AX$ the space of homogeneous polynomials of degree $1$. We still denote by $\varphi$ its linear extension to 
$AX\otimes AX$ given by 
\begin{equation}
	\varphi(P,Q)=\sum_{x,y\in X} \scal{P}{x}\scal{Q}{y} \varphi(x,y)
\end{equation}
and $\pshuffle{\varphi}$ the extension of the mapping of Definition \mref{Defstarproduct} by linearity\footnote{We recall that $AX$ (resp. $\ncp{A}{X}$) admits $X$ (resp. $X^*$) as linear basis, therefore $AX\otimes AX$ (resp. $\ncp{A}{X}\otimes \ncp{A}{X}$) is free with basis $X\times X$ (resp. $X^*\times X^*$) or more precisely, the image family $(x\otimes y)_{x,y\in X}$ (resp. $(u\otimes v)_{u,v\in X^*}$).} to $\ncp{A}{X}\otimes \ncp{A}{X}$. Then 
$\pshuffle{\varphi}$ becomes a law of algebra (with $1_{X^*}$ as unit) on $\ncp{A}{X}$.

\subsection{Extending quasi-stuffle relations}\label{ss:stuffl}

The following elementary result can be found in any complex analysis textbook. It is freely used throughout this section.  \begin{lemma}
For any integers $s,r\geq 1$, for any complex numbers $a$, $b\not=a$ :
\begin{eqnarray}
\forall x\in\C\setminus\{a,b\}, 
\dfrac1{(x-a)^s(x-b)^r}=\sum_{k=1}^s\dfrac{a_k}{(x-a)^k}+\sum_{k=1}^r\dfrac{b_k}{(x-b)^k}
\end{eqnarray}
where, for all $k\in\{1,\ldots,s\}$, 
$\displaystyle a_k=\binom{s+r-k-1}{r-1}\dfrac{(-1)^{s-k}}{(a-b)^{s+r-k}}$\\
and, for all $k\in\{1,\ldots,r\}$, 
$\displaystyle b_k=\binom{s+r-k-1}{s-1}\dfrac{(-1)^{r-k}}{(b-a)^{s+r-k}}$.

\end{lemma}

Let $\bt=(t_1,\ldots,t_r)$ be a set of parameters$^{\ref{ft:parameters}}$, $\bS=(s_1,\ldots,s_r)$ a composition, $\bxi=(\xi_1,\ldots,\xi_r)\in\C^r$. We define, for $N\in\N_{>0}$,
\begin{eqnarray}
M_{\bS,\bxi,\bt}^N&=&\sum_{N\geqslant n_1>\ldots>n_r>0}\prod_{i=1}^r\dfrac{\xi_i^{n_i}}{(n_i-t_i)^{s_i}}.
\end{eqnarray}
and $M^N_{(),(),()}=1$.\\[2pt]
Of course, it is a truncated series of $\zeta(\bS;\bt;\bxi)$.
\begin{proposition}\label{pro:M}
For every composition $\bS$, tuple $\bxi$ of complex numbers, tuple $\bt$ of parameters all of the same length $l\in\N$, and for every composition $\br$, tuple $\brho$ of complex numbers, tuple $\bt'$ of parameters also of the same length $k\in\N$, one has
\begin{eqnarray}
\forall N\in\N,\quad
M_{\bS,\bxi,\bt}^N\,M_{\br,\brho,\bt'}^N
=
M_{(\bS,\bxi,\bt)\luffle(\br,\brho,\bt')}^N.
\end{eqnarray}
\end{proposition}
\Proof
If $l=0$ or $k=0$, that is immediate.\\
Let $l\in\N^*$, $k\in\N^*$ and $\bS=(s_1,\ldots,s_l)$ and $\br=(r_1,\ldots,r_k)$ two compositions,   $\bxi=(\xi_1,\ldots,\xi_l)\in\C^l$,  $\brho=(\rho_1,\ldots,\rho_k)\in\C^k$, and $\bt=(t_1,\ldots,t_l)$, $\bt'=(t_1',\ldots,t_k')$ two sets of parameters and
put  $\bS_2=(s_2,\ldots,s_l)$,  $\br_2=(r_2,\ldots,r_k)$, $\bxi_2=(\xi_2,\ldots,\xi_l)$,  $\brho_2=(\rho_2,\ldots,\rho_k)$, $\bt_2=(t_2,\ldots,t_l)$ and ${\bt'}_2=({t'}_2,\ldots,{t'}_k)$, 
\begin{itemize}
\item[$\bullet$] If $t'_1=t_1$, 
\begin{eqnarray}
M_{\bS,\bxi,\bt}^N\,M_{\br,\brho,\bt}^N
&=&\sum_{N\geqslant n_1,N\geqslant {n'_1}}
\dfrac{\xi_1^{n_1}}{{(n_1-t_1)}^{s_1}}
\,M_{\bS',\bxi',\bt_2}^{n_1}\,
\dfrac{{\rho}_1^{{n'}_1}}{({n'}_1-{t}_1)^{r_1}}
\,M_{\br',\brho',{\bt}_2}^{{n'}_1}
\end{eqnarray}
Classically, we decompose the sum $\ds\sum_{N\geqslant n_1,N\geqslant >{n'_1}}$ into three sums corresponding to the simplices $n_1>n'_1;\ n_1'>n_1$ and $n_1=n'_1$ and get   
\begin{eqnarray}
M_{\bS,\bxi,\bt}^N\,M_{\br,\brho,\bt'}^N
&=&\sum_{N\geqslant n_1}
\dfrac{\xi_1^{n_1}}{(n_1-t_1)^{s_1}}
\,M_{\bS_2,\bxi_2,\bt_2}^{n_1}
\,M_{\br,\brho,\bt}^{n_1}
+\sum_{N\geqslant {n'}_1}
\dfrac{{\rho}_1^{{n'}_1}}{({n'}_1-{t'}_1)^{r_1}}
\,M_{\bS,\bxi,\bt}^{{n}_1}
\,M_{\br_2,\brho_2,{\bt}_2}^{{n'}_1}\cr
&&+\sum_{N\geqslant m}
\dfrac{(\xi_1\rho_1)^{m}}{(m-t_1)^{s_1}}\dfrac1{(m-{t}_1)^{r_1}}
\,M_{\bS_2,\bxi_2,\bt_2}^{m}
\,M_{\br_2,\brho_2,\bt'_2}^{m}
\end{eqnarray}
so that, 
\begin{eqnarray}
\forall N\in\N,\quad
M_{\bS,\bxi,\bt}^N\,M_{\br,\brho,\bt}^N
=
M_{(\bS,\bxi,\bt)\luffle(\br,\brho,\bt')}^N.
\end{eqnarray}
\item[$\bullet$] In the same way, when $t_1\neq t'_1$ 
\vskip-6mm
\begin{eqnarray}
M_{\bS,\bxi,\bt}^N\,M_{\br,\brho,\bt'}^N
&=&\sum_{N\geqslant n_1}
\dfrac{\xi_1^{n_1}}{(n_1-t_1)^{s_1}}
\,M_{\bS_2,\bxi_2,\bt_2}^{n_1}
\,M_{\br,\brho,{\bt'}}^{n_1}
+\sum_{N\geqslant {n'}_1}
\dfrac{{\rho}_1^{{n'}_1}}{({n'}_1-{t'}_1)^{r_1}}
\,M_{\bS,\bxi,{\bt}}^{{n'}_1}
\,M_{\br_2,\brho_2,{\bt'}_2}^{{n'}_1}\cr
&&+\sum_{N\geqslant m}
(\xi_1\rho_1)^{m}
\dfrac1{(m-t_1)^{s_1}}\dfrac1{(m-{t'}_1)^{r_1}}
\,M_{\bS_2,\bxi_2,\bt_2}^{m}
\,M_{\br_2,\brho_2,\bt_2}^{m}\cr
&=&\sum_{N\geqslant n_1}
\dfrac{\xi_1^{n_1}}{(n_1-t_1)^{s_1}}
\,M_{\bS_2,\bxi_2,\bt_2}^{n_1}
\,M_{\br,\brho,{\bt'}}^{n_1}
+\sum_{N\geqslant {n'}_1}
\dfrac{{\rho}_1^{{n'}_1}}{({n'}_1-{t'}_1)^{r_1}}
\,M_{\bS,\bxi,{\bt}}^{{n'}_1}
\,M_{\br_2,\brho_2,{\bt}_2}^{{n'}_1}\cr
&&+\sum_{N\geqslant m}
\Big\lbrack
\sum_{k=1}^{s_1}
\binom{s_1+r_1-k-1}{r_1-1}\dfrac{(-1)^{s_1-k}}{(t_1-t_1')^{s_1+r_1-k}}\dfrac{(\xi_1\rho_1)^{m}}{(m-{t'}_1)^k}\cr
&&+\sum_{k=1}^{r_1}
\binom{s_1+r_1-k-1}{s_1-1}\dfrac{(-1)^{r_1-k}}{(t_1-{t'}_1)^{s_1+r_1-k}}
\Big\rbrack
\,M_{\bS_2,\bxi_2,\bt}^{m}
\,M_{\br_2,\brho_2,\bt'}^{m}
\end{eqnarray}
\end{itemize}
so
\begin{eqnarray}
\forall N\in\N,\quad
M_{\bS,\bxi,\bt}^N\,M_{\br,\brho,\bt'}^N
=
M_{(\bS,\bxi,\bt)\luffle(\br,\brho,\bt')}^N.
\end{eqnarray}
\cqfd
\begin{remark}
Let $r$ a integer, $\chi=(\chi_1,\ldots,\chi_r)$ a tuple of multiplicative characters\footnote{Endomorphisms of the semigroup $(\C,\times)$.} and ($\bS,\bxi,\bt$ being as above) let us define
\begin{eqnarray}
M_{\bS,\bxi,\bt}^N(\chi)&=&\sum_{N\geqslant n_1>\ldots>n_r>0}\prod_{i=1}^r\dfrac{\chi_i^{n_i}(\xi_i)}{(n_i-t_i)^{s_i}}.
\end{eqnarray}
The same proof shows that, for any $(\bS,\bxi)\in\Z_{>0}^l\times\C^l$ and 
$(\br,\brho)\in\Z_{>0}^k\times\C^k$, for any l-tuple $\bt$ and $k-$tuple $\bt'$ of parameters$^{\ref{ft:parameters}}$,
\begin{eqnarray}
\forall N\in\N,\quad
M_{\bS,\bxi,\bt}^N(\chi)\,M_{\br,\brho,\bt'}^N(\chi)
=
M_{(\bS,\bxi,\bt)\luffle(\br,\brho,\bt')}^N(\chi).
\end{eqnarray}
\end{remark}
This result allows to deduce some product relations on the different multi-z\^eta functions :
\begin{theorem}\label{Th:stuffle}
Let  $\bS=(s_1,\ldots,s_l)$ and $\br=(r_1,\ldots,r_k)$ two compositions,  $\bxi=(\xi_1,\ldots,\xi_l)$ a $l$-tuple, $\brho=(\rho_1,\ldots,\rho_k)$ a $k$-tuple of complex numbers of which the first composant has a modulus strictly less than $1$, $\bt=(t_1,\ldots,t_s)$ and $\bt'=(t_1',\ldots,t_k')$ two tuples of parameters not in $\N_{>0}$ and $N\in \N$ 
\begin{itemize}
\item[(i)] For the coloured polyz\^eta function :
\begin{eqnarray}
\zeta(\bS,\bxi)\zeta(\bS',\bxi')&=&\zeta\left((\bS,\bxi)\duffle(\bS',\bxi')\right)
\end{eqnarray}
\item[(ii)] For the truncated Hurwitz polyz\^eta  function :
\begin{eqnarray}
\zeta_{_{N}}(\bS,\bt)\zeta_{_{N}}(\bS',\bt')&=&\zeta_{_{N}}\left((\bS,\bt)\huffle(\bS',\bt')\right)
\end{eqnarray}
\item[(iii)] In particular, for the monocentered polyz\^eta function :
\begin{eqnarray}
\zeta\left(\bS,(t,\ldots,t)\right)\zeta\left(\bS',(t,\ldots,t)\right)&=&\zeta\left(\left(\bS,(t,\ldots,t)\right)\stuffle\left(\bS',(t,\ldots,t)\right)\right)
\end{eqnarray}
where $t$ is a parameter s.t. $t\not\in\N_{>0}$.
\item[(iv)] For the Polylerch generalized function :
\begin{eqnarray}
\zeta(\bS,\bt,\bxi)\zeta(\bS',\bt',\bxi')&=&\zeta\left((\bS,\bt,\bxi)\luffle(\bS',\bt',\bxi')\right)
\end{eqnarray}
\end{itemize}
\end{theorem}
\Proof
(ii) comes directly from Proposition \ref{pro:M} because $\zeta_{_{N}}(\bS,\bt)=M^N_{\bS,(1,\ldots,1),\bt}$;
for (i), (iii) and (iv), apply Proposition \ref{pro:M} with, respectively, the functions
$$M^N_{\bS,\bxi,(0,\ldots,0)},\quad M^N_{\bS,(1,\ldots,1),(t,\ldots,t)}\quad\mbox{ and } M^N_{\bS,\bxi,\bt}$$ 
and take both sides of the equality to the limit as $N$ grows to infinity.
\cqfd
\begin{remark}
We cannot use this method for the Hurwitz polyz\^etas because in the decomposition, some divergent terms (which have $s_1=1$ !) appear: for example, for $t\neq t'$,
\begin{eqnarray}
(y_2,z_t)\huffle(y_3,z_{t'})
&=&(y_2y_3,z_tz_{t'})+(y_3y_2,z_{t'}z_t)
+\sum_{n=0}^1\binom{2+n}{2}\dfrac{(-1)^n}{(t-t')^{3+n}}(y_{2-n},z_t)\cr
&&+\sum_{n=0}^2\binom{1+n}{1}\dfrac{(-1)^n}{(t-t')^{2+n}}(y_{3-n},z_{t'})
\cr
&=&(y_2y_3,z_tz_{t'})+(y_3y_2,z_{t'}z_t)
+\dfrac{1}{(t-t')^3}(y_2,z_t)-\dfrac3{(t-t')^4}(y_1,z_t)\cr
&&+\dfrac{1}{(t-t')^2}(y_3,z_{t'})-\dfrac2{(t-t')^3}(y_1,z_{t'})+\dfrac{3}{(t-t')^4}(y_1,z_{t'})
\end{eqnarray}
Separately, 
the terms $-\dfrac3{(t-t')^4}(y_1,z_t)$ and $\dfrac{3}{(t-t')^4}(y_1,z_{t'})$, corresponding respectively to
$\dfrac{-3}{(t-t')^4}\dfrac1{n-t}$ and $\dfrac3{(t-t')^4}\dfrac1{n-t'}$ give a divergent series although all other terms correspond to convergent series. Of course, the sum of the two
\begin{eqnarray}
\dfrac3{(t-t')^4}\left(-\dfrac1{n-t}+\dfrac1{n-t'}\right)
=\dfrac3{(t-t')^4}\left(\dfrac{t'-t}{(n-t)(n-t')}\right)
\end{eqnarray}
is a term of a convergence series, but the series is not a Hurwitz Polyz\^eta.
\end{remark}

\section{Radford's theorem for the AC-stuffle.}
In this subsection, $A$ is supposed to be a ring with unit; when we need it to be commutative or to contain the set of rational numbers, we will state it explicitly.

Let $<$ be a total ordering on the alphabet $X$, and $\Lyn(X)$ denote the family of Lyndon words \cite{reutenauer} constructed from $X^*$ w.r.t. this ordering. 
We will prove that the largest framework in which Radford's theorem holds true \cite{radford} is when $\varphi$ is commutative (and associative).
\subsection{Computing $\varphi$-shuffle expressions using shuffles}
In this subsection $A$ is a ring with unit and $\varphi : AX\otimes AX\rightarrow AX$ an associative law.\\
We can express the result of the $\varphi$-shuffle product thanks to the shuffle product (and some terms of lower degree).
First we observe what happens with the product of two words :
\begin{lemma}\label{lem:passageShuffleDeuxMots}
For $u,v\in X^*$, there exists $(C_{u,v}^w)_{|w|<|u|+|v|}\in A^{(\N)}$ such that :
$$u\,\pshuffle{\varphi}v=u\shuffle v+\sum_{|w|<|u|+|v|}C_{u,v}^ww.$$
\end{lemma}
\Proof Omitted.
\cqfd

Now, because the  Lyndon words are candidates to be a transcendental basis, we see what happens when they are $\varphi$-shuffled.
\begin{definition}\label{multidegree_notation}
Let $\star:\ncp{A}{X}\times \ncp{A}{X}\mapsto\ncp{A}{X}$ be an associative law with unit and $\X=Lyn(X)$. For any $\alpha\in \N^{(\X)}$ and $\{l_1,\cdots ,l_r\}\supset supp(\alpha)$ in strict decreasing order (i.e. $l_1>\cdots >l_r$), we set 
\begin{equation}\label{multprod}
\X^{\star\alpha}= l_1^{\star\alpha_1}\star\cdots \star l_r^{\star\alpha_r}
\end{equation}  
\end{definition}
One easily checks easily that the product (\ref{multprod}) does not depend on the choice of $\{l_1,\cdots ,l_r\}\supset supp(\alpha)$. We will also need the following parameter (which will turn out to be the length of the dominant terms in the product) 
\begin{equation}
||\alpha||=\sum\limits_{l\in\Lyn(X)}\alpha(l)|l|\ .
\end{equation}

\begin{lemma}\label{lem:triangulaire}
If $\pshuffle{\varphi}$ is associative,
$$\displaystyle\forall \alpha\in\N^{(\Lyn(X))}, \exists(C_\beta^\alpha)_{\beta}\in A^{(\N^{(\Lyn(X))})}/ \,\X^{\pshuffle{\varphi}\alpha}=\X^{\shuffle\alpha}+\sum_{\stackrel{\beta\in\N^{(\Lyn(X))}}{||\beta||<||\alpha||}}C_\beta^\alpha\X^{\shuffle\beta}.$$
\end{lemma}
\Proof Omitted
\cqfd
\subsection{Radford's theorem in $\varphi$-shuffle algebras}
\begin{lemma}\label{generator}
If $\pshuffle{\varphi}$ is associative,
\begin{equation}\label{geneq}
\forall p\in\N^*, 
\span\left((\X^{\pshuffle{\varphi}\alpha})_{\alpha\in\N^{\left(\Lyn(X)\right)}, ||\alpha||<p}\right)
=\span\left((\X^{\shuffle\alpha})_{\alpha\in\N^{\left(\Lyn(X)\right)}, ||\alpha||<p}\right)
.
\end{equation}

\end{lemma}
\Proof
Lemma \ref{lem:triangulaire} give $\forall p\in\N^*, 
\span\left((\X^{\pshuffle{\varphi}\alpha})_{\alpha\in\N^{\left(\Lyn(X)\right)}, ||\alpha||<p}\right)
\subset\span\left((\X^{\shuffle\alpha})_{\alpha\in\N^{\left(\Lyn(X)\right)}, ||\alpha||<p}\right)
$.\\
We just have to prove, for any $p\in\N^*$, the property $\calP(p)$ :
\begin{eqnarray}
\span\left((\X^{\shuffle\alpha})_{\alpha\in\N^{\left(\Lyn(X)\right)}, ||\alpha||<p}\right)\subset\span\left((\X^{\pshuffle{\varphi}\alpha})_{\alpha\in\N^{\left(\Lyn(X)\right)}, ||\alpha||<p}\right)
\end{eqnarray}
\begin{itemize}
\item[$\bullet$] It is true for $p=1$.
\item[$\bullet$]  Assume $\calP(p)$ true for an integer $p$.\\
Let $\alpha\in\N^{\left(\Lyn(X)\right)}$ such that $||\alpha||<p+1$.\\
We can find $(C_\beta^\alpha)_{\beta}\in A^{(\N^{(\Lyn(X))})}$ such that 
$\X^{\pshuffle{\varphi}\alpha}=\X^{\shuffle\alpha}+\sum_{\stackrel{\beta\in\N^{(\Lyn(X))}}{||\beta||<||\alpha||}}C_\beta^\alpha\X^{\shuffle\beta}$, so
$\X^{\shuffle\alpha}=\X^{\pshuffle{\varphi}\alpha}-\sum_{\stackrel{\beta\in\N^{(\Lyn(X))}}{||\beta||<||\alpha||}}C_\beta^\alpha\X^{\shuffle\beta}$.\\
But every term of the sum is of the form $C_\beta^\alpha\X^{\shuffle\beta}$ with $\beta\in\N^{(\Lyn(X))}$ and $||\beta||<||\alpha||<p+1$ so $||\beta||<p$.\\
Consequently, they are in $\span\left((\X^{\pshuffle{\varphi}\alpha})_{\alpha\in\N^{\left(\Lyn(X)\right)}, ||\alpha||<p}\right)$, and so is the sum. By the induction hypothesis, 
the sum is in $\span\left((\X^{\shuffle{\varphi}\alpha})_{\alpha\in\N^{\left(\Lyn(X)\right)}, ||\alpha||<p}\right)$, therefore $\X^{\shuffle\alpha}\in\span\left((\X^{\shuffle{\varphi}\alpha})_{\alpha\in\N^{\left(\Lyn(X)\right)}, ||\alpha||<p+1}\right)$.
\end{itemize}
\cqfd
\begin{theorem}\label{th:Lyn-linbasis} 
Let $A$ be a commutative ring (with unit) such that $\Q\subset A$\footnote{Precisely, $\N^+.1_A\subset A^\times$} and $\pshuffle{\varphi} : \ncp{A}{X}\otimes \ncp{A}{X}\rightarrow \ncp{A}{X}$ is associative.\\
If $X$ is totally ordered by $<$, then $\left(\X^{\pshuffle{\varphi}\alpha}\right)_{\alpha\in\N^{(\Lyn(X))}}$ is a linear basis of $\ncp{A}{X}$.
\end{theorem}
\Proof
Since this family is a generating family by lemma \ref{generator}, only freeness remains to be proven.\\
Let $\displaystyle \sum_{\alpha\in J}\beta_{\alpha}\X^{\pshuffle{\varphi}\alpha}=0$ be a null linear combination of $(\X^{\pshuffle{\varphi}\alpha})_{\alpha\in\N^{\left(\Lyn(X)\right)}}$, with $J$ a nonempty finite subset of $\N^{\left(\Lyn(X)\right)}$.
Thanks to lemma \ref{lem:triangulaire}, for any $\alpha\in J$, we can find a finite family $B_{\alpha}\subset \N^{(\Lyn(X))}$ and $\displaystyle (C_\beta^\alpha)_{\beta\in B_{\alpha}}\in A^{B_{\alpha}}$ such that 
$$\X^{\pshuffle{\varphi}\alpha}=\X^{\shuffle\alpha}+\sum_{\stackrel{\beta\in B_{\alpha}}{||\beta||<||\alpha||}}C_\beta^\alpha\X^{\shuffle\beta}.$$
Set $\displaystyle B=J\cup\left(\bigcup_{\alpha\in J}B_{\alpha}\right)$; $B$ is a finite set.
Then $(\X^{\pshuffle{\varphi}\alpha})_{\alpha\in J}$ is a triangular family for $|.|$ with respect to the family ${\mathcal F}=(\X^{\pshuffle{\varphi}\alpha})_{\alpha\in B}$ in the vector space $\span({\mathcal F})$, which is of finite dimension. But ${\mathcal F}$ is a basis, so $(\X^{\pshuffle{\varphi}\alpha})_{\alpha\in J}$ is free and 
$\forall \alpha\in J, \beta_{\alpha}=0$.
\cqfd
\begin{corollary}
Under the same hypotheses, if in addition $\pshuffle{\varphi}$ is commutative in $A$ 
then 
\begin{itemize}
\item[i)] The algebra $\mathcal{A}=(\ncp{A}{X},\pshuffle{\varphi},1_{X^*})$ is a polynomial algebra. 
\item[ii)] $Lyn(X)$ is a transcendence basis of $\mathcal{A}$.
\end{itemize}
\end{corollary}
\begin{remark}
It is necessary to suppose $\Q\subset A$ as, in case $\varphi\equiv 0$, 
one has 
\begin{equation}
\forall n\in\N_{>0},	a^n=\frac{1}{n!}(a^{\shuffle n})
\end{equation}
\end{remark}
\Proof
\begin{itemize}
\item[i)] Immediate result.
\item[ii)] Comes from proposition \ref{th:Lyn-linbasis} and theorem \ref{th:prop_law}, which proves in an elementary (so independent) way that 
 the commutativity of $\varphi$ is equivalent to  the commutativity of $\pshuffle{\varphi}$.
\end{itemize}
\cqfd
\subsection{Bialgebra structure}
\begin{definition}
A law $\star$ defined over $\ncp{A}{X}$ is a dual law (or dualizable) if there exists   
a linear mapping 
$\Delta_{\star}:\ncp{A}{X}\rightarrow\ncp{A}{X}\otimes\ncp{A}{X}$ such 
\begin{eqnarray}
\forall(u,v,w)\in X^*\times X^*\times X^*, \quad
\scal{u\star v}{w}=\scal{u\otimes v}{\Delta_{\star}(w)}^{\otimes 2}\ .
\end{eqnarray}
In this case, $\Delta_*$ will be called the comultiplication dual to $\star$.
\end{definition}
\begin{theorem}\label{bialgebra} If $A$ is a commutative ring (with unit), if $\Q\subset A$, and if in addition the product $\shuffle_{\varphi} : \ncp{A}{X}\otimes \ncp{A}{X}\rightarrow \ncp{A}{X}$ is an associative and commutative law on $\ncp{A}{X}$, then
the algebra $(\ncp{A}{X},\pshuffle{\varphi},1_{X^*})$ can be endowed with 
the comultiplication $\Delta_{\conc}$ dual to the concatenation
\begin{equation}
	\Delta_{\conc}(w)=\sum_{uv=w} u\otimes v
\end{equation}
and the ``constant term'' character $\epsilon(P)=\scal{P}{1_{X^*}}$.\\
With this setting  
\begin{equation}\label{b1}
	\mathcal{B}_{\varphi}=(\ncp{A}{X},\pshuffle{\varphi}, 1_{X^*}, \Delta_{\conc},\epsilon) 
\end{equation}
is a bialgebra \footnote{Commutative and, when $|X|\geq 2$, noncocommutative.}.
\end{theorem}
\begin{remark} 
Let, classically, $\Delta_{\conc}^+$ be defined by 
$$
\forall w\in X^*, \Delta_{\conc}^+(w)=\sum_{uv=w\atop u,v\not=1} u\otimes v
$$
We remark that $\Delta_{\conc}^+$ is coassociative and locally nilpotent, i.e.
$$\forall w\in X^*,\exists n\in\N^*/  \left(\Delta_{\conc}^+\right)^n(w)=0.$$
Thus the bialgebra \mref{b1} is, in fact, a Hopf Algebra.
\end{remark}
\Proof
It is a classical combinatorial verification, done in \cite{DM12}.\\ 
The following identity remains to be proven: 
 \begin{eqnarray}
 \forall (w_1,w_2)\in X^*, \Delta_{\conc}(w_1\pshuffle{\varphi}w_2)=\Delta_{\conc}(w_1)\Delta_{\conc}(w_2)
 \end{eqnarray}
which can be done by a (lengthy) induction or by duality. 
\cqfd

\section{Conditions for AC-shuffle and dualizability}

\subsection{Commutative and associative conditions}

We have obtained an extended version of Radford's theorem and other properties with conditions stated w.r.t. $\pshuffle{\varphi}$, we will see in this subsection that these conditions can be set uniquely in terms of properties of $\varphi$ itself. 
\begin{definition}
For $P\in\ncp{A}{X}$, we note $supp(P)$ the support of $P$ and  
$\displaystyle \deg(P)=\mbox{max}\{|l|,l\in supp(P)\}$
\end{definition}
\begin{lemma}
Let $A$ be a unitary commutative ring, $X$ be an alphabet and $\varphi: X\times X\rightarrow \ncp{A}{X}$ is an arbitrary mapping.\\
Then, 
\begin{eqnarray}
\forall(u,v)\in(X^*)^2, \deg(u\pshuffle{\varphi}v)\pp|u|+|v|
\end{eqnarray}
\end{lemma}
\Proof
If $|u|=0$ or $|v|=0$, then $u\pshuffle{\varphi}v$ is one of $\{u,v\}$ so its length is $|u|+|v|$.\\
Let $X^+$ be the set of nonempty words. We prove 
$\forall(u,v)\in(X^+)^2, \deg(u\pshuffle{\varphi}v)=|u|+|v|$  
by induction on $|u|+|v|$.\\
For any letters $a$ and $b$, $a\pshuffle{\varphi}b=ab+ba+\varphi(a,b)1_{A^*}$
so $\deg(a\pshuffle{\varphi}b)=2=|a|+|b|$.\\
One assumes the property true for all words $u,v \in X^+$ such that $|u|+|v|=n$, where $n$ is an integer. Let $u$ and $v$ be now two words of $X^+$ such that  $|u|+|v|=n+1$.\\
There exist $x, y$ in $X$, $u'$, $v'$ in $X^*$ such that $u=xu'$, $v=yv'$ (because $(u, v)\in(X^+)^2$).
Then $|u|+|v'|=|u|+|v|-1\pp n$, so $\deg(y(u\pshuffle{\varphi}v'))\pp n+1$. Also 
$|u'|+|v|=|u|-1+|v|\pp n$ so $\deg(x(u'\pshuffle{\varphi}v))\pp n+1$, and 
$|u'|+|v'|=|u|-1+|v|-1\pp n$ so $\deg(\varphi(x,y)u'\pshuffle{\varphi}v')\pp n+1$.
Hence, $\deg(u\pshuffle{\varphi}v)=n+1$ : the induction is proved. 
\cqfd
\begin{theorem}\label{th:prop_law}
In the context of definition \ref{Defstarproduct},
\begin{itemize}
\item[(i)] The law $\pshuffle{\varphi}$ is commutative if and only if the extension $\varphi : AX\otimes AX\rightarrow AX$ is commutative. 
\item[(ii)] The law $\pshuffle{\varphi}$ is associative if and only if the extension $\varphi : AX\otimes AX\rightarrow AX$ is associative.
\end{itemize}
\end{theorem}
\Proof We give an elementary proof.
\begin{itemize}
\item[(i)] [$\pshuffle{\varphi}$ commutative $\Longrightarrow$ $\varphi$ commutative]\\  
Let us suppose\quad
$
\forall (u,v)\in {(X^*)}^2,\quad u\pshuffle{\varphi} v=v\pshuffle{\varphi} u.
$.\\
In particular,  $\forall (x,y)\in {(X^*)}^2, x\pshuffle{\varphi} y=x\pshuffle{\varphi} y.$
But, for any $(x,y)\in X^2$ , 
\begin{eqnarray}
x\pshuffle{\varphi} y=xy+yx+\varphi(x,y)&\mbox{and}&y\pshuffle{\varphi} x=yx+xy+\varphi(y,x).
\end{eqnarray}
and so
$(\forall x,y\in X)(\varphi(x,y)=\varphi(y,x))$.

[$\varphi$ commutative $\Longrightarrow$ $\pshuffle{\varphi}$ commutative]\\
Now let us suppose $\varphi$ is commutative then let us prove by recurrence on $|uv|$ that $\pshuffle{\varphi}$ is commutative :
\begin{itemize}
\item The previous equivalence proves that the recurrence holds for $|u|=|v|=1$.
\item Suppose the recurrence holds for any $u,v\in X^*$ such that $2\le|uv|\le n$ and $|u|,|v|\neq1$.\\
Let $u=xu'$ and $v=yv'$ with $x,y\in X$ and $u',v'\in X^*$. Then,
\begin{eqnarray}
u\pshuffle{\varphi} v
&=&x(u'\pshuffle{\varphi} yv)+y(xu'\pshuffle{\varphi} v)+\varphi(x,y)(u'\pshuffle{\varphi} v')\cr
&=&x(yv\pshuffle{\varphi} u)+y(v'\pshuffle{\varphi} xu')+\varphi(y,x)(v'\pshuffle{\varphi} u')\qquad(\mbox{by the induction hypothesis})\cr
&=&v\pshuffle{\varphi} u.
\end{eqnarray}
\end{itemize}

\item[(ii)] [$\pshuffle{\varphi}$ associative $\Longrightarrow$ $\varphi$ associative] Let us suppose
\begin{eqnarray}
\forall u,v,w\in X^*,&&(u\pshuffle{\varphi} v)\pshuffle{\varphi} w=u\pshuffle{\varphi}(v\pshuffle{\varphi} w).
\end{eqnarray}
Then, for any $x,y,z\in X$ , one has
\begin{eqnarray}
(x\pshuffle{\varphi} y)\pshuffle{\varphi} z=x\pshuffle{\varphi}(y\pshuffle{\varphi} z).
\end{eqnarray}
But
\begin{eqnarray}
(x\pshuffle{\varphi} y)\pshuffle{\varphi} z&=&(xy+yx+\varphi(x,y))\pshuffle{\varphi} z\\
&=&xy\pshuffle{\varphi} z+yx\pshuffle{\varphi} z+\varphi(x,y)\pshuffle{\varphi} z\cr
&=&x(y\pshuffle{\varphi} z)+z(xy\pshuffle{\varphi} 1)+\varphi(x,z)y+y(x\pshuffle{\varphi} z)+z(yx\pshuffle{\varphi} 1)+\varphi(y,z)x\cr 
&+& \varphi(x,y)z+z\varphi(x,y)+\varphi(\varphi(x,y),z)\cr
&=&x(yz+zy+\varphi(y,z))+zxy+\varphi(x,z)y+y(xz+zx+\varphi(x,z))+zyx\cr
&+&\varphi(y,z)x+\varphi(x,y)z+z\varphi(x,y)+\varphi(\varphi(x,y),z)
\cr\cr
x\pshuffle{\varphi} (y\pshuffle{\varphi} z)&=& x\pshuffle{\varphi} (yz+zy+\varphi(y,z))\\
&=&x\pshuffle{\varphi} yz+x\pshuffle{\varphi} zy+x\pshuffle{\varphi} \varphi(y,z)\cr
&=& x(1\pshuffle{\varphi} yz)+ y(x\pshuffle{\varphi} z)+ \varphi(x,y)z+x(1\pshuffle{\varphi} zy)+ z(x\pshuffle{\varphi} y)+ \varphi(x,z)y\cr
&=& x\varphi(y,z)+\varphi(y,z)x+\varphi(x,\varphi(y,z))\cr	
&=& xyz+y(xz+zx+\varphi(x,z))+ \varphi(x,y)z+xzy+ z(xy+yx+\varphi(x,y))+\varphi(x,z)y\cr
&+& x\varphi(y,z)+\varphi(y,z)x+\varphi(x,\varphi(y,z)) .
\end{eqnarray}
One can then deduce that
\begin{eqnarray}
(\forall x,y,z\in X)(x\pshuffle{\varphi} (y\pshuffle{\varphi} z)=(x\pshuffle{\varphi} y)\pshuffle{\varphi} z)
&\iff&
(\forall x,y,z\in X)(\varphi(x,\varphi(y,z))=\varphi(\varphi(x,y),z)).\cr
\end{eqnarray}
[$\varphi$ associative $\Longrightarrow$ $\pshuffle{\varphi}$ associative] 
Now if $\varphi$ is associative then let us prove by induction on $|u|+|v|+|w|$ that $\pshuffle{\varphi}$ is associative :
\begin{itemize}
\item The previous equivalence proves that the induction holds for $|u|=|v|=|w|=1$.
\item Suppose the recurrence holds for any $u,v\in X^*$ such that $3\le|u|+|v|+|w|\le n$ and $|u|,|v|,|w|\neq1$.
\item Let $u=xu',v=yv'$ and $w=zw'$ with $x,y,z\in X$ and $u',v',w'\in X^*$. Then,
\begin{eqnarray}
&&u\pshuffle{\varphi}(v\pshuffle{\varphi} w)\cr
&&=u\pshuffle{\varphi}\left(y(v'\pshuffle{\varphi} w)+z(v\pshuffle{\varphi} w')+\varphi(y,z)(v'\pshuffle{\varphi} w')\right)\cr
&&=x(u'\pshuffle{\varphi}y(v'\pshuffle{\varphi} w))+y(u\pshuffle{\varphi}(v'\pshuffle{\varphi} w))+\varphi(x,y)(u'\pshuffle{\varphi}(v'\pshuffle{\varphi} w))\cr
&&+x(u'\pshuffle{\varphi}z(v\pshuffle{\varphi} w'))+z(u\pshuffle{\varphi}(v\pshuffle{\varphi} w'))+\varphi(x,z)(u'\pshuffle{\varphi}(v\pshuffle{\varphi} w'))\cr
&&+x(u'\pshuffle{\varphi}\varphi(y,z)(v'\pshuffle{\varphi} w'))+\varphi(y,z)(u\pshuffle{\varphi}(v'\pshuffle{\varphi} w'))+\varphi(x,\varphi(y,z))u'\pshuffle{\varphi}(v'\pshuffle{\varphi} w')\cr
&&=x(u'\pshuffle{\varphi}(v\pshuffle{\varphi} w))\cr
&&+y(u\pshuffle{\varphi}(v'\pshuffle{\varphi} w))+\varphi(x,y)(u'\pshuffle{\varphi}(v'\pshuffle{\varphi} w))\cr
&&+z(u\pshuffle{\varphi}(v\pshuffle{\varphi} w'))+\varphi(x,z)(u'\pshuffle{\varphi}(v\pshuffle{\varphi} w'))\cr
&&+\varphi(y,z)(u\pshuffle{\varphi}(v'\pshuffle{\varphi} w'))+\varphi(x,\varphi(y,z))u'\pshuffle{\varphi}(v'\pshuffle{\varphi} w')\\
\mbox{and}&&\cr
&&(u\pshuffle{\varphi}v)\pshuffle{\varphi} w\cr
&&=(x(u'\pshuffle{\varphi}v)+y(u\pshuffle{\varphi}v')+\varphi(x,y)(u'\pshuffle{\varphi}v'))\pshuffle{\varphi} w)\cr
&&=x((u'\pshuffle{\varphi}v)\pshuffle{\varphi} w)+z(x(u'\pshuffle{\varphi}v)\pshuffle{\varphi} w')+\varphi(x,z)((u'\pshuffle{\varphi}v)\pshuffle{\varphi} w')\cr
&&+y((u\pshuffle{\varphi}v')\pshuffle{\varphi} w)+z(y(u\pshuffle{\varphi}v')\pshuffle{\varphi} w')+\varphi(y,z)((u\pshuffle{\varphi}v')\pshuffle{\varphi} w')\cr
&&+\varphi(x,y)((u'\pshuffle{\varphi}v')\pshuffle{\varphi} w)+z(\varphi(x,y)(u'\pshuffle{\varphi}v')\pshuffle{\varphi} w')+\varphi(\varphi(x,y),z)((u'\pshuffle{\varphi}v')\pshuffle{\varphi} w')\cr
&&=x((u'\pshuffle{\varphi}v)\pshuffle{\varphi} w)+\varphi(x,z)((u'\pshuffle{\varphi}v)\pshuffle{\varphi} w')\cr
&&+y((u\pshuffle{\varphi}v')\pshuffle{\varphi} w)+\varphi(y,z)((u\pshuffle{\varphi}v')\pshuffle{\varphi} w')\cr
&&+\varphi(x,y)((u'\pshuffle{\varphi}v')\pshuffle{\varphi} w)+\varphi(\varphi(x,y),z)((u'\pshuffle{\varphi}v')\pshuffle{\varphi} w')\cr
&&+z(u\pshuffle{\varphi}v)\pshuffle{\varphi} w')
\end{eqnarray}
Indeed, thanks to the induction hypothesis and the commutativity of $\varphi$, $u\pshuffle{\varphi}(v\pshuffle{\varphi} w)$ and $(u\pshuffle{\varphi}v)\pshuffle{\varphi} w$ are equal.
\end{itemize}
\end{itemize}
\cqfd
\subsection{Dualizability conditions}

\begin{proposition}
We call $\gamma_{x,y}^z:=\scal{\varphi(x,y)}{z}$ the structure constants of $\varphi$ (w.r.t. the basis $X$).\\ 
The product $\pshuffle{\varphi}$ is a dual law if and only if $(\gamma_{x,y}^z)_{x,y,z\in X}$ is dualizable in the following sense
\begin{eqnarray}
	(\forall z\in X)(\#\{(x,y)\in X^2|\gamma_{x,y}^z\not=0\}<+\infty)\ .
\end{eqnarray}
\end{proposition}
\Proof
[$\pshuffle{\varphi}$ dual law $\Longrightarrow$ $\gamma_{x,y}^z$ dualizable]. Let $\Delta$ be the dual of $\pshuffle{\varphi}$, that is, for all $u,v,w\in X^*$ 
\begin{equation}
	\scal{u\pshuffle{\varphi} v}{w}=\scal{u\otimes v}{\Delta(w)}^{\otimes 2}\ .
\end{equation}
For all $z\in X$, one must have $\Delta(z)=\sum_{i=1}^n \alpha_i u_i\otimes v_i$. On the other hand, for all $x,y\in X$, one has $(x\pshuffle{\varphi} y)-(xy+yx)=\varphi(x,y)$. Hence 
\begin{eqnarray}
\gamma_{x,y}^z&=&\scal{\varphi(x,y)}{z}=\scal{(x\pshuffle{\varphi} y)-(xy+yx)}{z}=
\scal{(x\pshuffle{\varphi} y)}{z}-\scal{(xy+yx)}{z}\cr
&&=\scal{(x\otimes y)}{\Delta(z)}=
\scal{(x\otimes y)}{\sum_{i=1}^n \alpha_i u_i\otimes v_i}\ .
\end{eqnarray}
We can deduce from the preceding argument that 
$$
\gamma_{x,y}^z\not=0\Longrightarrow \left(x\in \cup_{i=1}^n Alph(u_i)\textrm{ and }y\in \cup_{i=1}^n Alph(v_i)\right)
$$
which proves the point. 

\smallskip
[$\gamma_{x,y}^z$ dualizable $\Longrightarrow$ $\pshuffle{\varphi}$ dual law]). This is, combinatorially speaking, the most interesting point. We first define a comultiplication $\Delta$ on $\ncp{A}{X}$ by transposing the structure constants of $\pshuffle{\varphi}$ by 
\begin{equation}\label{letters_Delta}
	\Delta(z):=z\otimes 1+1\otimes z+\sum_{x,y\in X}\gamma_{x,y}^z x\otimes y
\end{equation}
and, as the sum is finite (see however the comment after this theorem), this quantity belongs to $\ncp{A}{X}\otimes \ncp{A}{X}$. One then has a linear mapping $\Delta : AX\rightarrow \ncp{A}{X}\otimes \ncp{A}{X}$ which is extended, by universal property, into a morphism of algebras $	\Delta : \ncp{A}{X}\rightarrow \ncp{A}{X}\otimes \ncp{A}{X}$. Explicitely, for all $w=z_1z_2\cdots z_n$, one has 
\begin{equation}\label{ext_Delta}
\Delta(z_1z_2\cdots z_n)=\Delta(z_1)\Delta(z_2)\cdots \Delta(z_n)\ .
\end{equation}
Now, we prove that the dual law of the latter coproduct is exactly $\pshuffle{\varphi}$.\\ 
First remark : by \mref{letters_Delta} and \mref{ext_Delta}, one has 
\begin{equation}
\Delta(w)=w\otimes 1+1\otimes w+\sum_{u,v\in X^+} \beta(u,v) u\otimes v     
\end{equation}
the last sum being finitely supported. This shows by duality that 
\begin{equation}
u\pshuffle{\Delta} 1=1\pshuffle{\Delta} u=u	
\end{equation}
(here, $\pshuffle{\Delta}$ stands for the dual law of $\Delta$). Moreover
\begin{eqnarray}
au\pshuffle{\Delta} bv &=& \sum_{w\in X^*}\scal{au\pshuffle{\Delta} bv}{w}w\cr
&=&\sum_{w\in X^*}\scal{au\otimes bv}{\Delta(w)}w\cr
&=&\scal{au\otimes bv}{1\otimes 1} 1+\sum_{w\in X^+}\scal{au\otimes bv}{\Delta(w)}w\cr
&=&
\sum_{x\in X\, ; \, m\in X^*}\scal{au\otimes bv}{\Delta(xm)}\ xm\cr
&=&\sum_{x\in X\, ; \, m\in X^*}\scal{au\otimes bv}{\Delta(x)\Delta(m)}\ xm\cr
&=&\sum_{x\in X\, ; \, m\in X^*}\scal{au\otimes bv}
{\Big(x\otimes 1+1\otimes x+\sum_{y,z\in X}\scal{\Delta(x)}{y\otimes z}\, y\otimes z\Big)\Delta(m)}\ xm\cr
&=&\sum_{x\in X\, ; \, m\in X^*}\scal{au\otimes bv}
{(x\otimes 1)\Delta(m)}\ xm+
\sum_{x\in X\, ; \, m\in X^*}\scal{au\otimes bv}
{(1\otimes x)\Delta(m)}\ xm+\cr
&&\sum_{x\in X\, ; \, m\in X^*}\scal{au\otimes bv}
{\sum_{y,z\in X}\scal{\Delta(x)}{y\otimes z}\, y\otimes z)\Delta(m)}\ xm\cr
&&\sum_{m\in X^*}\scal{au\otimes bv}{(a\otimes 1)\Delta(m)}\ am+
\sum_{m\in X^*}\scal{au\otimes bv}{(1\otimes b)\Delta(m)}\ bm+\cr
&&\sum_{x\in X\, ; \, m\in X^*}\scal{au\otimes bv}
{\scal{\Delta(x)}{a\otimes b}\, a\otimes b)\Delta(m)}\ xm\cr
&&\sum_{m\in X^*}\scal{u\otimes bv}{\Delta(m)}\ am+
\sum_{m\in X^*}\scal{au\otimes v}{\Delta(m)}\ bm+\cr
&&\sum_{x\in X\, ; \, m\in X^*}\scal{\Delta(x)}{a\otimes b}\scal{u\otimes v}
{\Delta(m)}\ xm\cr
&=&a\sum_{m\in X^*}\scal{u\otimes bv}{\Delta(m)}\ m+
b\sum_{m\in X^*}\scal{au\otimes v}{\Delta(m)}\ m+\cr
&&\sum_{m\in X^*}\Big(\sum_{x\in X}\scal{\Delta(x)}{a\otimes b}\ x\Big)\scal{u\otimes v}
{\Delta(m)}\ m\cr
\cr
&=&a(u\pshuffle{\Delta} bv)+b(au\pshuffle{\Delta} v)+\varphi(a,b)(u\pshuffle{\Delta} v)
\end{eqnarray}
This proves that the dual law $\pshuffle{\Delta}$ equals $\pshuffle{\varphi}$ and we are done. 
\cqfd

\subsection{The Hopf-Hurwitz algebra}

In section \ref{HP}, we provided the law on indices followed by the product of Formal Hurwitz polyz\^etas, we now prove that the law $\varphi$ associated with it is associative. The ``centres'' will be taken from a subfield $k$ of $\C$ and the set of coefficients $A$ is a $k$-CAAU.
\begin{proposition}
\begin{enumerate}
\item[i)] The law $\varphi : AN\otimes AN\rightarrow AN$ associated to $\circ$ is defined, on the basis $N$, by the multiplication table $\mathcal{T}_{Formal\ Hurwitz}$
\begin{eqnarray}
&\mbox{ if }  t= t';&\ \varphi((y_i,z_t),(y_j,z_{t'}))=(y_{i+j},z_{t})\cr
&\mbox{ if }  t\not= t';&\ \varphi((y_i,z_t),(y_j,z_{t'}))=\cr
&&\ \displaystyle\sum_{n=0}^{i-1}{j-1+n\choose j-1} 
\frac{(-1)^n}{(t-t')^{j+n}}\,(y_{i-n},z_t)
+\displaystyle\sum_{n=0}^{j-1} {i-1+n\choose i-1} 
\frac{(-1)^n}{(t'-t)^{i+n}}\,(y_{j-n},z_{t'})\cr
\end{eqnarray}  
\item[ii)] The product $\huffle$ is associative, commutative and unital, making $(\ncp{A}{N},\circ,1_N)$ into a $A$-CAAU.
\end{enumerate}
\end{proposition}
\Proof
i) Let first $j : kN\rightarrow k(X)$ be the linear mapping defined by $j((y_i,z_t))=\frac{1}{(X-t)^i}$. In fact, as the $\{\frac{1}{(X-t)^i}\}$ are linearly independent, $j$ is into. On the other hand, $j$ is a morphism of $k$-AAU due to the fact that the multiplication table is identical. Hence $\varphi$ is a law of $\C$-AA on $AN\simeq A\otimes_k kN$.
%
ii) Is a consequence of the general theorems.      
\cqfd  

Now, we have the following bialgebra
\begin{eqnarray}
\calH_{Formal\ Hurwitz}&=&(\ncp{A}{N},{\ministuffle_{\varphi}},1_{N^*},\Delta_{\tt conc},\epsilon)\label{Hurwitz-Hopf}
\end{eqnarray}    
which is a Hopf algebra. Note that $\ministuffle_{\varphi}$ is not dualisable which means that the adjoint 
\begin{equation}
\Delta_{\ministuffle_{\varphi}} : N^*\rightarrow \ncs{A}{N^*\otimes N^*}
\end{equation}  
does not have its image in $\ncp{A}{N}\otimes \ncp{A}{N}$. See next paragraph for tools and proofs. 
\begin{corollary}
The product $\luffle$ is associative, commutative and unital, making $(\ncp{A}{N},\luffle,1_N)$ into a $A$-CAAU.
\end{corollary}
\Proof
It comes that the product $\luffle$ is a direct product of the products $\circ$ and $\smuffle$.
\cqfd

\section{Conclusion}

We have been able to give a useful extended version of Radford's theorem.\\
Let us observe that :
\begin{itemize}
\item[$\bullet$] For the shuffle product, $\varphi_{\shuffle}\equiv0$, so the shuffle $\shuffle$ is associative, commutative and dualizable.
\item[$\bullet$] The stuffle product over an alphabet indexed by $\N$ is associative and commutative because  $\varphi_{\stuffle}(x_i,x_j)=x_{i+j}$ is so; moreover it is dualizable.
\item[$\bullet$] The muffle product  over an alphabet indexed by $\C$ is associative and commutative because  $\varphi_{\smuffle}(x_i,x_j)=x_{i\times j}$ is so; it is not dualizable because 
for all $n\in\N_{>0}$, $x_1=\varphi_{\smuffle}(x_{1/n},x_n)$.\\
However, there are multiplicative subsemigroups $S$ of $\C$ such that $\varphi$ restricted to the alphabet $(x_i)_{i\in S}$ is dualizable.
\item[$\bullet$] The duffle product  over an alphabet indexed by $\N^*\times\C^*$   is associative and commutative because $\varphi_{\duffle}\big((y_i,x_k),(y_j,x_l)\big)=(y_{i+j},x_{k\times l})$ is associative and commutative; it is not dualizable either (for the same reason).\\
But we can do the same remark as the muffle about the possibility to restrict the alphabet so that $\varphi$ becomes dualizable.
\item[$\bullet$] The Formal Polyz\^eta product is associative, commutative but, for all $j\in\N_{>0}$, with $t\neq t'$,
\begin{eqnarray}
\varphi_{_{\circ}}\left((y_2,z_t),(y_j,z_{t'})\right)
&=&\sum_{n=0}^1\binom{j-1+n}{j-1}\dfrac{(-1)^n}{(t-t')^{j+n}}(y_{2-n},z_t)\cr
&&+\sum_{n=0}^{j-1}\binom{1+n}{1}\dfrac{(-1)^n}{(t'-t)^{2+n}}(y_{j-n},z_{t'})
\end{eqnarray}
and then $\displaystyle\qquad\forall j\in\N_{>0},\quad
\scal{\varphi_{_{\circ}}\left((y_2,z_t),(y_j,z_{t'})\right)}{(y_1,z_t)} = -\binom{j}{j-1}\dfrac1{(t-t')^{j+1}} \neq 0$\\
which implies that $\circ$ is not dualizable. 
\item[$\bullet$] The Lerch  product is associative, commutative and not dualizable (for the same reason as $\smuffle$).
\end{itemize}
So, if we work in the Riemann polyz\^eta algebra, in the coloured polyz\^eta algebra, 
or in the Generalized Lerch polyz\^eta algebra,
we can use a representation with the Lyndon set as a transcendence basis. 
Moreover, in the Riemann polyz\^eta algebra and the truncated Hurwitz polyz\^eta algebra can both be completed into Hopf algebras.

\end{document}